\documentclass[12pt]{article}

\usepackage{amssymb}    
\usepackage{amsthm}     

\newtheorem{theorem}{Theorem}[section]
\newtheorem{lemma}[theorem]{Lemma}
\newtheorem{proposition}[theorem]{Proposition}
\newtheorem{corollary}[theorem]{Corollary}

\newenvironment{prf}[1]{\trivlist
\item[\hskip
\labelsep{\it #1.\hspace*{.3em}}]}{
\endtrivlist}

\newtheorem{predefinition}[theorem]{Definition}
\newenvironment{definition}{\begin{predefinition}\rm}{\end{predefinition}}
\newtheorem{preremark}[theorem]{Remark}
\newenvironment{remark}{\begin{preremark}\rm}{\end{preremark}}
\newtheorem{prenotation}[theorem]{Notation}

\newtheorem{preexample}[theorem]{Example}
\newenvironment{example}{\begin{preexample}\rm}{\end{preexample}}
\newtheorem{preclaim}[theorem]{Claim}

\newtheorem{prequestion}[theorem]{Question}

 \makeatletter
\def\emppsubsection{\@startsection{subsection}{2}{\z@}{-3.25ex plus -1ex minus -.2ex}{-1em}{\bf}}

\makeatother

\newcommand \Y {{\cal Y}}
\newcommand \CO {{\cal O}}

\newcommand \PP {{\mathbb P}^1}

\newcommand \Aa {{\mathbb A}^1}
\newcommand \ZZ {{\mathbb Z}}
\newcommand \NN {{\mathbb N}}

\newcommand  \RR {{\mathbb R}}

\newcommand  \QQ {{\mathbb Q}}

\newcommand \dime {\mathop{\rm dim}}
\newcommand \Frac {\mathop{\rm Frac}}

\newcommand \gen {\mathop{\rm genus}}

\newcommand \Hom {\mathop{\rm Hom}}
\newcommand \Ind {\mathop{\rm Ind}}

\newcommand \mini {\mathop{\rm min}}

\newcommand \Spec {\mathop{\rm Spec}}

\newcommand \val {\mathop{\rm val}}

\title{Wildly ramified covers with large genus}
\author{Rachel J.\ Pries: Colorado State University
\footnote{The author was partially supported by NSF grant DMS-04-00461.}}
\date{}

\begin{document}
\maketitle

\begin{abstract}
\noindent
We study wildly ramified $G$-Galois covers $\phi:Y \to X$ branched at $B$
(defined over an algebraically closed field of characteristic $p$).
We show that curves $Y$ of arbitrarily high genus occur for such covers
even when $G$, $X$, $B$ and the inertia groups are fixed.
The proof relies on a Galois action on covers of germs of curves and formal patching. 
As a corollary, we prove that for any nontrivial quasi-$p$ group $G$ and 
for any sufficiently large integer $\sigma$ with $p \nmid \sigma$, 
there exists a $G$-Galois \'etale cover of the affine line with conductor $\sigma$
above the point $\infty$.

2000 Mathematical Subject Classification: 14H30, 14G32
\end{abstract}

\section{Introduction.}

Consider a $G$-Galois cover $\phi:Y \to X$ of smooth projective irreducible curves
defined over an algebraically closed field $k$ of characteristic $p$.
If $\phi$ is tamely ramified, the Riemann-Hurwitz formula implies that 
the genus of $Y$ is determined by $|G|$, $g_X$, the size of the branch locus $B$
and the orders of the inertia groups.  
This statement is no longer true when $\phi$ is wildly ramified.
Not only can wildly ramified covers usually be deformed without varying $X$ or the branch
locus $B$ of $\phi$, but they can be often be distinguished from each other by 
studying finer ramification invariants such as the conductor.
The genus of $Y$ now depends on these finer ramification invariants. 

This phenomenon is already apparent when one considers 
$\ZZ/p$-Galois covers $\phi:Y \to \PP_k$ of the projective line branched only over $\infty$.
Each of these can be given by an Artin-Schreier equation $y^p-y=f(x)$ 
where the degree $j$ of $f(x) \in k[x]$ is relatively prime to $p$.
The genus $g_Y$ of $Y$ equals $(p-1)(j-1)/2$ and thus can be arbitrarily large.  

The main point of this paper is that the same unboundedness phenomenon 
occurs for covers of any affine curve with any Galois group whose 
order is divisible by $p$.  
Namely, as in the case of $\ZZ/p$-covers of the affine line, 
the discrete invariants of covers of a fixed affine curve with fixed group 
can be arbitrarily large.
When $G$ is an abelian $p$-group, these results are a well-understood
consequence of class field theory.
Already the case of non-abelian $p$-groups does not seem to appear in the literature.
The complexity of the problem is related to the fact 
that the Galois group can be a simple group in which case the cover $\phi$
will not factor through a Galois cover of smaller degree.
For this reason, we use the technique of formal patching. 

Here are the main results.  
One discrete invariant for an \'etale cover of the affine line is its {\it conductor} 
(or last jump in the filtration of ramification groups in the upper numbering).
Abhyankar's Conjecture states that a $G$-Galois cover of the affine line exists 
if and only if $G$ is {\it quasi-$p$}, i.e.\ if and only if $G$ is generated by $p$-groups.

\vspace{.1in}
\noindent {\bf Corollary \ref{Ca1jump}.}
{\it If $G \not = 0$ is a quasi-$p$ group and $\sigma \in \NN$ ($p \nmid \sigma$)
is sufficiently large, then
there exists a $G$-Galois cover $\phi:Y \to \PP_k$ branched at only  
one point with conductor $\sigma$.}   

\smallskip 

More generally, suppose $X$ is a smooth projective irreducible $k$-curve 
and $B \subset X$ is a non-empty finite set of points.
Suppose $G$ is a finite quotient of $\pi_1(X-B)$ so that $p$ divides $|G|$.  
(These groups were classified by Raynaud \cite{Ra:ab} 
and Harbater \cite{Ha:ab} in their proof of Abhyankar's Conjecture.)
We show that curves $Y$ of arbitrarily high genus occur 
for $G$-Galois covers $\phi:Y \to X$ branched at $B$.

\vspace{.1in}
\noindent {\bf Corollary \ref{CbigG}.}
{\it There is an arithmetic progression ${\mathfrak p} \subset \NN$
so that for all $g \in {\mathfrak p}$,  
there exists a $G$-Galois cover $\phi:Y \to X$ branched only at $B$ with $\gen(Y)=g$.}

\smallskip

In addition, we give a lower bound for the proportion of natural numbers 
which occur as the genus of $Y$ for a $G$-Galois cover $\phi:Y \to X$ branched at $B$.  
This lower bound depends only on $G$ and $p$ and not on $X$ and $B$.
 
Hurwitz spaces for tamely ramified covers are well-understood by the work of 
\cite{W:thesis}.  The following corollary shows that the 
the structure of Hurwitz spaces for wildly ramified covers will be vastly different.  

\vspace{.1in}
\noindent {\bf Corollary \ref{Churwitz}}
{\it For any smooth projective irreducible $k$-curve $X$
and any non-empty finite set of points $B \subset X$
and any finite quotient $G$ of $\pi_1(X-B)$ so that $p$ divides $|G|$,   
a Hurwitz space for $G$-Galois covers $\phi:Y \to X$ branched at $B$ will have infinitely 
many components.}

\bigskip


The main results, all appearing in Section \ref{S3}, use formal patching to 
reduce to the question of deformations of a wildly ramified 
cover $\hat{\phi}$ of germs of curves.
We study these deformations by analyzing a Galois action on the $I$-Galois cover 
$\hat{\phi}$ in Section \ref{S1}.      
Here $I$ corresponds to the inertia group of $\phi$ at a ramification point 
and is of the form $P \rtimes_\iota \mu_m$ where $|P|=p^e$ and $p \nmid m$.
We use group theory and ramification theory to factor $\hat{\phi}$ 
as $\kappa \circ \overline{\phi} \circ \phi^{A}$.  The Galois group of $\phi^{A}$ is an 
elementary abelian $p$-group denoted $A$, the Galois group of $\overline{\phi}$ is $P/A$, 
and $\kappa$ is a Kummer $\mu_m$-cover.  
We define a Galois action on $\hat{\phi}$ 
which acts non-trivially only on the subcover $\phi^{A}$.
It causes the set of wildly ramified $I$-Galois covers of germs 
of curves which dominate $\overline{\phi} \circ \kappa$  
to form a principal homogeneous space under the action of an explicitly computable group.
The effect of this Galois action on the ramification filtration of a cover is investigated
in Proposition \ref{Ph1action}.

Another application of this Galois action can be found in \cite{Pr:deform}, 
where it is used to study deformations of wildly ramified
covers of germs of curves with control over the conductor.
The results in Section \ref{S1} appear in the generality needed for this application.

I would like to thank D. Harbater, along with J. Achter and I. Bouw,
for many suggestions which helped improve earlier drafts of this paper.  
I would also like to thank the participants of the conferences in Banff and in Leiden
for enlightening conversations.
 
\section{Galois action on wildly ramified covers.} \label{S1}

In this section, we investigate a Galois action on a 
wildly ramified cover of germs of curves and its effect on the conductor.    
The main result is Proposition \ref{Ph1action} which we use in Section \ref{S3} 
to control the change in 
the ramification filtration when we modify a Galois cover of germs of curves.

Let $k$ be an algebraically closed field of characteristic $p > 0$.
Consider a compatible system of roots of unity of $k$.  
Denote by $\zeta$ the chosen generator of $\mu_m$.
We start in Section \ref{S1.1} with background on the ramification filtrations 
and upper and lower jumps of a wildly ramified cover $\hat{\phi}$ of germs of curves. 
 

\subsection{Higher ramification groups.} \label{S1.1}

In this section, we define the filtration of higher ramification groups of a wildly ramified
Galois cover ${\hat{\phi}}:{\hat {Y}} \to U$ of germs of curves.
More precisely, we fix an irreducible $k$-scheme $\Omega$
and let $U=\underline{\Spec}(\CO_\Omega[[u]])$. 
Let $\xi$ be the point of height one of $U$ defined by the equation $u=0$.

Suppose ${\hat{\phi}}:{\hat {Y}} \to U$ is a Galois cover of normal connected germs 
of $\Omega$-curves which is wildly ramified at the closed point 
$\eta ={\hat{\phi}}^{-1}(\xi) \in {\hat {Y}}$.
By \cite[Lemma 2.1.4]{Pr:deg}, after an \'etale pullback of $\Omega$, 
the decomposition group and inertia group over the generic point of $\eta$ are the same
and so the Galois group of ${\hat{\phi}}$ is the same as the inertia group.  

Recall that the inertia group $I$ at the generic point of $\eta$ is of the form 
$P \rtimes_\iota \mu_m$ where $|P|=p^e$ for some $e > 0$ and $p \nmid m$, \cite{Se:lf}. 
Here $\iota$ denotes the automorphism of $P$ which determines the 
conjugation action of $\mu_m$ on $P$.

Associated to the cover ${\hat{\phi}}$, there are two filtrations of $I$, namely 
the filtration of higher ramification groups $I_{c'}$ in the lower numbering
and the filtration of higher ramification groups $I^c$ in the upper numbering. 
If $c' \in \NN$, then $I_{c'}$ is the normal subgroup of all $g \in I$ 
such that $g$ acts trivially on $\CO_{\hat {Y}}/\pi^{c'+1}$.
Here $\pi$ is a uniformizer of $\CO_{\hat {Y}}$ at the generic point of $\eta$. 
Equivalently, $I_{c'}=\{g \in I| \val(g(\pi)-\pi) \geq c'+1\}$. 
If $c' \in \RR^+$ and $c''=\lfloor c' \rfloor$, then $I_{c'}$ is defined to equal $I_{c''}$.   
Recall by Herbrand's formula \cite[IV, Section 3]{Se:lf}, 
that the filtration $I^c$ in the upper numbering is given by
$I^c=I_{c'}$ where $c'=\Psi(c)$ and $\Psi(c)=\int_0^c(I^0:I^t)dt$.  
Equivalently, $c=\varphi(c')$ where $\varphi(c')= \int_0^{c'}dt/(I_0:I_t)$.

If $I_{j} \not = I_{j+1}$ for some $j \in \NN^+$, we say that $j$ is a {\it lower jump}
of ${\hat{\phi}}$ at $\eta$.  The multiplicity of $j$ is the integer $\ell$ so that 
$I_j/I_{j+1} \simeq (\ZZ/p)^\ell$.  A rational number $c$ is an {\it upper jump}
of ${\hat{\phi}}$ at $\eta$ if $c=\varphi(j)$ for some lower jump $j$.   
We denote by $j_1, \ldots, j_e$ (resp.\ $\sigma_1, \ldots, \sigma_e$)
the set of lower (resp.\ upper) jumps of ${\hat{\phi}}$ at $\eta$
written in increasing order with multiplicity.  
Note that these are the positive breaks in the filtration of ramification groups in the lower 
(resp.\ upper) numbering.
By \cite[IV, Proposition 11]{Se:lf}, $p \nmid j_i$ for any lower jump $j_i$.
Herbrand's formula implies that $j_i-j_{i-1}=(\sigma_i-\sigma_{i-1})|I|/|I_{j_i}|$. 
In particular, $\sigma_i|I|/|I^{\sigma_i}| \in \NN$ and $\sigma_1=j_1/m$.

We call $\sigma=\sigma_e$ the {\it conductor} of 
${\hat{\phi}}$ at $\eta$; $\sigma$ is the largest $c \in \QQ$ such that 
inertia group $I^c$ is non-trivial in the filtration of higher ramification groups 
in the upper numbering.  (Note that this indexing is slightly different
than in \cite{Se:lf}, where if 
$x$ is a uniformizer at the branch point, then the 
ideal $(x^{\sigma+1})$ is the {\it conductor} of the 
extension of complete discrete valuation rings.) 

If $\phi:Y \to X$ is a $G$-Galois cover of projective curves branched at $B$, we briefly recall how
the genus $g_Y$ of $Y$ depends on the upper jumps of $\phi$ at each branch point.
Let $I_b$ denote the inertia group of $\phi$ at a point above $b \in B$
and let $mp^e$ denote its order; 
(for simplicity we drop the index $b$ from the variables $m$ and $e$).

\begin{lemma} \label{Lgenus} (Riemann-Hurwitz formula)
$$g_Y=1+|G|(g_X-1)+\sum_{b \in B} |G|{\rm deg}(R_b)/2|I_b|.$$
$${\rm deg}(R_b)=|I_b|-1+(p-1)m[\sigma_1 +p\sigma_2 + \ldots p^{e-1}\sigma_e].$$
\end{lemma}

\begin{proof}
The proof follows immediately from the Riemann-Hurwitz formula \cite[IV.2.4]{Har}, 
\cite[IV, Proposition 4]{Se:lf}, and Herbrand's formula \cite[IV.3]{Se:lf}.
\end{proof}

We will increase the conductor of ${\hat{\phi}}$ by modifying its $A$-Galois subcover while fixing 
its $I/A$-Galois quotient for a suitable choice of $A \subset I$.

\begin{lemma} \label{LhypA}  Suppose $I \simeq P \rtimes_\iota \mu_m$.
Suppose ${\hat{\phi}}:{\hat {Y}} \to U$ is an $I$-Galois cover with conductor $\sigma$.  
Then there exists $A \subset I^\sigma$ satisfying the following hypotheses: 
$A$ is central in $P$; $A$ is normal in $I$; 
$A$ is a nontrivial elementary abelian $p$-group; and $A$ is irreducible under the $\mu_m$-action.
\end{lemma}

\begin{proof}
By \cite[IV, Cor.\ 3 through Prop.\ 7]{Se:lf}, we see that $I^\sigma$ is elementary abelian.
By \cite[IV, Prop.\ 10]{Se:lf}, if $g \in I=I_1$ and 
$h \in I_{\Psi(\sigma)}$
then $ghg^{-1}h^{-1} \in I_{\Psi(\sigma) +1}={\rm id}$. 
Thus $I^\sigma$ is central in $P$.
The subgroup $I^\sigma$ is normal in $I$ by definition.  
Thus the conjugation action of $\mu_m$ stablizes $I^\sigma$. 
So any choice for $A$ among nontrivial subgroups of $I^\sigma$ 
stabilized and irreducible under the action of $\mu_m$ 
satisfies all the hypotheses. 
\end{proof}

Given an $I$-Galois cover ${\hat{\phi}}:{\hat {Y}} \to U$ with conductor $\sigma$, 
we fix $A \subset P$ satisfying the hypotheses of Lemma \ref{LhypA}.
Let $a$ be the positive integer such that $A \simeq (\ZZ/p)^a$. 
Let $A \rtimes_\iota \mu_m$ be the semi-direct product determined by the 
restriction of the conjugation action of $\mu_m$ on $P$.
We fix a set of generators 
$\{ \tau_i| \ 1 \leq i \leq a\}$ for $A$.  
Let $\overline{P}=P/A$ and $\overline{I}=I/A$.

The condition that $A$ is central in $P$ (resp.\ $A$ is normal in $I$)
is used to define a transitive action of $A$-Galois covers on $P$-Galois covers in 
Section \ref{S1.2} (resp.\ $I$-Galois covers in Section \ref{S1.3}).
The conditions that $A \simeq (\ZZ/p)^a$ and that $A$ is irreducible 
under the action of $\mu_m$ make it easy to describe these $A$-Galois 
covers with equations.  The condition that $A \subset I^\sigma$ is important for controling 
the conductor when performing this action which is necessary in Section \ref{S1.4}
for Proposition \ref{Ph1action}.   

\subsection{Definition of the Galois action.} \label{S1.2}

Suppose ${\hat{\phi}}:{\hat {Y}} \to U$ is an $I$-Galois cover 
and $A$ is a subgroup in the center of $P$.
This yields a factorization of ${\hat{\phi}}$ which we denote 
$${\hat {Y}} \stackrel{\phi^A}{\to} \overline{Y} \stackrel{\overline{\phi}}{\to} \overline{X} 
\stackrel{\kappa}{\to} U.$$
Here ${\phi^A}$ is $A$-Galois, $\overline{\phi}$ is $\overline{P}$-Galois,
and $\kappa$ is $\mu_m$-Galois.   
Let ${\hat{\phi}}^P:{\hat {Y}} \to \overline{X}$ denote the $P$-Galois subcover of ${\hat{\phi}}$. 
If $A$ is normal in $I$ and $\overline{I}=I/A$, then 
$\kappa \circ \overline{\phi}$ is an $\overline{I}$-Galois cover. 
Let $\overline{\eta}={\phi^A}(\eta)$ (resp.\ $\xi'=\kappa^{-1}(\xi)$) 
be the ramification point of $\overline{Y}$ (resp.\ $\overline{X}$).
So $$\eta \mapsto \overline{\eta} \mapsto \xi' \mapsto \xi.$$

Also, $\overline{X} \simeq \underline{\Spec}(\CO_\Omega)[[x]])$ for some $x$ such that 
$x^m=u$.  The generator $\zeta$ of $\mu_m$ acts as $\zeta(x)=\zeta_m^h(x)$ 
for some integer $h$ relatively prime to $m$.
Let $U'=U - \{\xi\}$ and $X'=\kappa^{-1}(U')=\overline{X}-\{\xi'\}$.

Consider the group $H_A=\Hom(\pi_1(X'), A)$. 
We suppress the choice of basepoint from the notation.
An element $\alpha \in H_{A}$ 
may be identified with the isomorphism class of 
an $A$-Galois cover of $\overline{X}$ branched only over the closed point $\xi'$.
We denote this cover by $\psi_\alpha:V \to \overline{X}$.

\begin{lemma} \label{Ltransitive1}
Let $H_{A}=\Hom(\pi_1(X'), A)$.  If $A$ is in the center of $P$, then
the fibre $H_{\overline{\phi}}$ of $\Hom(\pi_1(X'), P) \to 
\Hom(\pi_1(X'), \overline{P})$ over $\overline{\phi}$ is a 
principal homogeneous space for $H_{A}$.
In other words, $H_A$ acts simply transitively on the fibre $H_{\overline{\phi}}$.
\end{lemma}

\begin{proof}
If $\gamma \in H_{\overline{\phi}}$ and $\omega \in \pi_1(X')$, then 
we define the action of $\alpha \in H_{A}$ on 
$H_{\overline{\phi}}$ as follows:
$\alpha \gamma(\omega)=\alpha(\omega) \gamma(\omega) \in P$.
The action is well-defined since $\alpha \gamma \in \Hom(\pi_1(X'), P)$ and
the image of $\alpha \gamma (\omega)$ and $\gamma (\omega)$ in 
$\overline{P}$ are equal.
The action is transitive since if $\gamma, \gamma' \in H_{\overline{\phi}}$ 
then $\omega \to \gamma(\omega)\gamma'(\omega)^{-1}$ is an element of $H_A$ 
because $A$ is in the center of $P$.
The action is simple since $\alpha_1=\alpha_2$ if 
$\alpha_1(\omega)\gamma(\omega)= \alpha_2(\omega)\gamma(\omega)$ for all $\omega \in \pi_1(X')$.
\end{proof}

For more details about Lemma \ref{Ltransitive1} in the case that 
$A \simeq \ZZ/p$, see \cite[Section 4]{HS}.
We note that the fibre $H_{\overline{\phi}}$ may be 
identified with isomorphism classes
of (possibly disconnected) $P$-Galois covers ${\hat {Y}} \to \overline{X}$ of normal germs of curves
which are \'etale over $X'$ and
dominate the $\overline{P}$-Galois cover $\overline{\phi}:\overline{Y} \to \overline{X}$.

We now give a more explicit formulation of this action.  
Suppose $\alpha \in H_A$ corresponds to the isomorphism
class of an $A$-Galois cover
$\psi_\alpha: V \to \overline{X}$ and $\gamma \in \Hom(\pi_1(X'), P)$ 
corresponds to the isomorphism 
class of a $P$-Galois cover ${\hat{\phi}}^P_\gamma:{\hat {Y}} \to \overline{X}$.
How do $\psi_\alpha$ and $\hat{\phi}^P_\gamma$ determine the cover corresponding 
to $\alpha \gamma$, which we denote $\hat{\phi}^P_{\alpha \gamma}$?  
Consider the normalized fibre products $Z={\hat {Y}} \tilde{\times}_{\overline{X}} V$
and $\overline{Z}=\overline{Y} \tilde{\times}_{\overline{X}} V$. 

\[\begin{array}{ccc}
Z&\stackrel{1 \times A}{\longrightarrow}&{\hat {Y}}\\
\downarrow{P}&&\downarrow{{\hat{\phi}}_{\gamma}^P}\\
V&\stackrel{\psi_{\alpha}}{\longrightarrow}&\overline{X}
\end{array}\]

Now $({\hat{\phi}}^P_\gamma \tilde{\times} \psi_\alpha): Z \rightarrow \overline{X}$ is the 
$P \times A$-Galois cover corresponding to the element
$(\gamma, \alpha) \in \Hom(\pi_1(X'), P \times A)$. 
Let $A' \subset P \times A$ be the normal subgroup generated by 
$(\tau_i, \tau_i^{-1})$ for $1 \leq i \leq a$. 
Note that $A' \simeq A$ and that $P \simeq (P \times A)/A'$.   
We denote by ${\hat{\phi}}^P_{\alpha \gamma}: W \rightarrow \overline{X}$ the
$P$-Galois quotient $({\hat{\phi}}^P_\gamma \tilde{\times} \psi_\alpha)^{A'}$
corresponding to the fixed field $W$ of $A'$ in $Z$.

\[\begin{array}{ccc}
Z&\stackrel{1 \times A}{\longrightarrow}&{\hat {Y}}\\
\downarrow{A'}&&\downarrow{{\hat{\phi}}_{\gamma}^A}\\
W&\stackrel{(A\times A)/A'}{\longrightarrow}&\overline{Y}
\end{array}\]

\begin{lemma} \label{Lfibre}
If $A$ is in the center of $P$, then the $P$-Galois cover corresponding to 
$\alpha \gamma \in \Hom(\pi_1(X'), P)$ is isomorphic to 
${\hat{\phi}}^P_{\alpha \gamma}:W \to \overline{X}$.
\end{lemma}

\begin{proof}
If $\omega \in \pi_1(X')$ then $\alpha \gamma(\omega)=
\alpha(\omega) \gamma(\omega) \in P$.
On the other hand, 
$(\gamma, \alpha)(\omega)=(\gamma(\omega), \alpha(\omega)) \in P \times A$.
Since $(\gamma(\omega), \alpha(\omega)) \equiv 
(\gamma(\omega)\alpha(\omega), 1)$ modulo
$A'=\langle (\tau_i, \tau_i^{-1}) \rangle$,
the two covers give the same element of $\Hom(\pi_1(X'), P)$.
\end{proof}

\subsection{Invariance under the $\mu_m$-action.} \label{S1.3}

In this section, we consider the invariance of these covers under the $\mu_m$-Galois
action.  Suppose $A$ is central in $P$ and normal in $I$.
Since $\iota$ restricts to an automorphism of $A$, we see that
$\iota$ acts naturally on $H_A=\Hom(\pi_1(X'), A)$.   
Denote by $H_A^\iota$ the subgroup of $H_A$ fixed by $\iota$.
In other words, the elements of $H_A^\iota$ correspond to $A$-Galois covers
$\psi:V \to \overline{X}$ branched only over the closed point $\xi'$ for which the composition 
$\kappa \circ \psi:V \to U$ is an $(A \rtimes_\iota \mu_m)$-Galois cover.

Suppose ${\hat{\phi}}: {\hat {Y}} \to U$ is an $I$-Galois cover 
where $I \simeq P \rtimes_\iota \mu_m$.  
Let $\gamma \in \Hom(\pi_1(X'), P)$ correspond to the isomorphism class of 
the $P$-Galois subcover ${\hat{\phi}}^P$ of ${\hat{\phi}}$.
If $\alpha \in H_A$, recall that ${\hat{\phi}}^P_{\alpha \gamma}$ is  
the cover corresponding to $\alpha \gamma \in \Hom(\pi_1(X'), P)$.

\begin{lemma}   
Suppose $A$ is central in $P$ and normal in $I$.  Let $\alpha \in H_A$.
The cover $\kappa \circ 
({\hat{\phi}}^P_{\alpha \gamma}):W \to U$ is an $I$-Galois cover if and only if 
$\alpha \in H_A^\iota$. 
\end{lemma}

\begin{proof}
Let $(P\times A) \rtimes_\iota \mu_m$ be the semi-direct product 
for which the conjugation action of $\zeta$ on 
the subgroups $(0, A)$ and $(A,0)$ in $P\times A$ is the same.  
By the hypotheses on $A$, 
the subgroup $A'=(A,A^{-1})$ is normal in $(P \times A) \rtimes_\iota \mu_m$.
So by Lemma \ref{Lfibre}, the cover 
$\kappa \circ ({\hat{\phi}}^P_{\alpha \gamma}):W \to U$ is $I$-Galois if and only if 
$Z \to U$ is $(P\times A) \rtimes_\iota \mu_m$-Galois.
One direction is now immediate: if $\alpha \in H_A^\iota$, 
then the cover $Z \to U$ is $(P \times A) \rtimes_\iota \mu_m$-Galois 
so $\kappa \circ ({\hat{\phi}}^P_{\alpha \gamma}):W \to U$ is $I$-Galois.

Conversely, suppose $\kappa \circ {\hat{\phi}}^P_{\alpha \gamma}:W \to U$ is $I$-Galois.
Since ${\hat {Y}} \to U$ and $W \to U$ are $I$-Galois 
the action of $\mu_m$ extends to an automorphism of ${\hat {Y}}$ and of $W$,   
which reduce to the same automorphism of $\overline{Y}$. 
Since ${\hat {Y}}$ is the fixed field of $Z$ under $(0,A) \subset P\times A$
and $W$ is the fixed field of $Z$ under $A'$, we see that 
the action of $\mu_m$ extends to an automorphism of $Z$.
Thus $Z \to U$ is $(P \times A)\rtimes_\iota \mu_m$-Galois.  
Recall that $V$ is the quotient of $Z$ by the normal subgroup $(P,0)$ of
$(P \times A)\rtimes_\iota \mu_m$.  
So the quotient $V \to U$ is $(A \rtimes_\iota \mu_m)$-Galois which implies 
$\alpha \in H_A^\iota$.  
\end{proof}

\subsection{The irreducible elementary abelian case.} \label{S1.35}

In this section, 
we consider the case that $A \subset P$ is a non-trivial elementary abelian $p$-group
$(\ZZ/p)^a$ which is irreducible under the action of $\mu_m$ on $P$. 
In other words, $A$ satisfies all the hypotheses of Lemma \ref{LhypA}
except that $A \subset I^\sigma$.
In this case, the results in Sections \ref{S1.2} and \ref{S1.3} can be made more explicit.  
Let $\delta_{ij}=1$ if $i=j$ and $\delta_{ij}=0$ otherwise.

First, if $A$ is a non-trivial elementary abelian $p$-group, 
then the cover $\psi_\alpha$ corresponding to $\alpha \in H_A$
is determined by its $\langle \tau_i \rangle$-Galois quotients
which are given by equations $v_i^{p}-v_i=r_{\alpha_i}$ with 
$r_{\alpha_i} \in \CO_\Omega[[x]][x^{-1}]$.  The generator $\tau_i \in A$ corresponds
to an automorphism $\tau_{\psi,i}$ of $V$ given by 
$\tau_{\psi,i}(v_j)=v_j +\delta_{ij}$.
Given another such cover $\psi'$ with equations 
$v_i^{p}-v_i=r'_i$, 
the group operation for $H_A$ corresponds to adding the Laurent series
in the Artin-Schreier equations.  This yields the cover
given by the equations $v^p_i-v_i=r_{\alpha_i}+r'_i$.

Next, if $A$ is a non-trivial elementary abelian $p$-group, 
then Lemma \ref{Lfibre} allows us to view the action of $H_A$ on 
$H_{\overline{\phi}}$ on the ring level.  The cover ${\hat{\phi}}^P_\gamma$ 
corresponding to $\gamma \in \Hom(\pi_1(X'), P)$ 
is determined from its $\overline{\phi}$-quotient by an $A$-Galois cover 
$\phi^A_\gamma:{\hat {Y}} \to \overline{Y}$.  This cover $\phi^A_\gamma$
is determined by its $\langle \tau_i \rangle$-Galois quotients
which are given by Artin-Schreier equations $y_i^{p}-y_i=r_{\gamma_i}$.
Here $r_{\gamma_i} \in \overline{K}$ where $\overline{K}$ is the 
fraction field of the complete local ring of $\overline{Y}$ at the closed point $\overline{\eta}$.  
The Galois action is given by $\tau_{{\hat{\phi}}, i}(y_j)=y_j +\delta_{ij}$.
Note that $w_i=y_i +v_i$ is invariant under 
$\langle (\tau_i, \tau_i^{-1}) \rangle$.   
Thus the cover ${\hat{\phi}}^P_{\alpha \gamma}:W \to \overline{X}$ 
is determined from $\overline{\phi}$ by the $A$-Galois cover
$W \to \overline{Y}$ which is determined by the equations 
$w_i^{p}-w_i= r_{\gamma_i} +r_{\alpha_i}$.   
This gives an even more explicit description of the Galois action.

Furthermore, suppose $A \subset P$ is a non-trivial elementary abelian $p$-group
so that $\mu_m$ acts irreducibly on $A$.  If $\psi_\alpha$ is the $A$-Galois cover corresponding 
to $\alpha \in H_A^\iota$, then $\psi_\alpha$ is determined by any one of its 
$\langle \tau_i \rangle$-Galois quotients along with the cover $\kappa$.  
Also the filtration of higher ramification groups for $\alpha \in H_A^\iota$
can have only one jump which occurs with full multiplicity $a$.  
In addition, the conductor of $\psi_\alpha$ is an integer which must satisfy a certain 
congruence condition.

\begin{lemma} \label{Lcong}
\begin{description}
\item{i)}
Suppose ${\hat{\phi}}: {\hat {Y}} \to U$ is an $I \simeq P \rtimes_\iota \mu_m$-Galois cover with 
conductor $\sigma$ and last lower jump $j_e$.  
Suppose $A \subset I^\sigma$ satisfies the hypotheses of Lemma \ref{LhypA}.
Suppose $\alpha \in H_A^\iota$ has conductor $s$.  
Then $s \equiv j_e/ |\overline{P}| \bmod m$.   
\item{ii)}
Associated to the $\mu_m$-Galois cover $\kappa: \overline{X} \to U$ and the group 
$A \rtimes_\iota \mu_m$, there is a unique integer $s_\iota$ 
(such that $1 \leq s_\iota \leq m$) having the following property:
if $\alpha \in H_A^\iota$ has conductor $s$ then $s \equiv s_\iota \bmod m$.
\end{description}
\end{lemma}

 
\begin{proof}
\begin{description}
\item{i)}
Suppose $|\overline{P}|=p^d$.
Let $y$ be a primitive element for the $A$-Galois extension of the function 
field $\overline{K}$ of $\overline{Y}$. 
The valuation of $y$ in the function field $\hat{K}$ of ${\hat {Y}}$ is $-j_e$ where $j_e$ 
is the last lower jump of ${\hat{\phi}}$. 
The equation for the $A$-Galois cover ${\phi^A}: {\hat {Y}} \to \overline{Y}$ 
is of the form $f(y)=r_{\hat{\phi}}$ where 
$f(y) \in \CO_\Omega[y]$ is a relative Eisenstein polynomial of degree $p^a$ 
and $r_{\hat{\phi}} \in \overline{K}$.
From the equation, the valuation of $r_{\hat{\phi}}$ in $\hat{K}$ is $p^a \val(y)=-p^a j_e$.
It follows that the valuation of $r_{\hat{\phi}}$ in $\overline{K}$ is $-j_e$.
Also $x$ has valuation $p^d$ in $\overline{K}$.  
This implies that there is a unit $u_1$ of $\CO_{\overline{K}}$
so that $r_{\hat{\phi}}=x^{-j_e/p^d}u_1$.  
Likewise, the equation for the $A$-Galois cover $\psi:V \to \overline{X}$ 
is of the form $f(v)=r_\alpha$ where $r_\alpha=x^{-s}u_2$ for some unit $u_2$ 
of $\CO_\Omega[[x]]$.   
The generator $\zeta$ of $\mu_m$ acts by $\zeta(x)=\zeta_m^hx$.  
So $\zeta(v)=\zeta_m^{-hs}v$ and $\zeta(y)=\zeta_m^{-hj_e/p^d}y$.
Since $y+v$ is in the function field of $W$, the action of $\zeta$ 
on $y$ and $v$ must be compatible.
Thus $s \equiv j_e/ |\overline{P}| \bmod m$.

\item{ii)} It is sufficient to show that the conductors $s_i$ of any two 
covers $\alpha, \alpha_1 \in H_A^\iota$ are congruent modulo $m$.  This follows 
directly from part (i) (taking ${\hat{\phi}}=\kappa \circ \alpha_1$, 
$P=A$, and $j_e=s_1$).  
\end{description}
\end{proof}

\subsection{Effect of the Galois action on the conductor.} \label{S1.4}

At this point, we determine the effect of the Galois action
on the conductor of the cover.  
Consider an $I$-Galois cover ${\hat{\phi}}:{\hat {Y}} \to U$ with conductor $\sigma$
and, more generally, upper jumps 
$\sigma_1, \ldots, \sigma_e$ in the ramification groups $I^c_{\hat{\phi}}$ 
in the upper numbering; (here $\sigma=\sigma_e)$.  
Consider $A \subset I^{\sigma}_{\hat{\phi}}$ satisfying 
the hypotheses of Lemma \ref{LhypA} and choose $\alpha \in H^\iota_A$.
By definition, $\alpha$ corresponds to an $A$-Galois cover $\psi_\alpha:V \to \overline{X}$.  
Let $s \in \NN^+$ be the conductor of $\psi_\alpha$ and note $p \nmid s$.
Since $A$ is irreducible under the $\mu_m$-action, 
$\psi_\alpha$ has ramification filtration $I^c_\alpha =A$ for $0 \leq c \leq s$ and 
$I^c_\alpha =0$ for $c > s$. 
The action of $\alpha$ takes ${\hat{\phi}}$ to another $I$-Galois cover
which we denote by ${\hat{\phi}}^{\alpha}:W \to U$.
Recall that ${\hat{\phi}}^\alpha$ is the $I$-Galois cover 
$\kappa \circ ({\hat{\phi}}^P_{\alpha \gamma})$ where 
$\gamma \in \Hom(\pi_1(X'), P)$ corresponds to the $P$-Galois subcover 
${\hat{\phi}}^P$.

The following result will be crucial in order to
modify a Galois cover of germs of curves with control over the change in 
the ramification filtration.

\begin{proposition} \label{Ph1action}
If the cover ${\hat{\phi}}^\alpha$ is connected,  
then it has conductor $\max\{s/m, \sigma\}$.
More generally, the ramification filtration of ${\hat{\phi}}^{\alpha}$
is $I_{{\hat{\phi}}^\alpha}^c=I_{{\hat{\phi}}}^c$ for $0 \leq c \leq \sigma$,
$I_{{\hat{\phi}}^\alpha}^c=A$ for $\sigma < c \leq \max\{s/m, \sigma\}$,
and $I_{{\hat{\phi}}^\alpha}^c=0$ for $\max\{s/m, \sigma\} < c$.
\end{proposition}

\begin{proof} 
The $P$-Galois subcover 
${\hat{\phi}}^P_{\alpha \gamma}=({\hat{\phi}}^{\alpha})^P: W \to \overline{X}$ dominates
$\overline{\phi}:\overline{Y} \to \overline{X}$.  So after intersecting with 
$\overline{P}$, the ramification filtrations
for ${\hat{\phi}}^P$ and ${\hat{\phi}}^P_{\alpha \gamma}$ are equal, i.e.\
$I_{{\hat{\phi}}^{\alpha}}^c \cap \overline{P}=I_{{\hat{\phi}}}^c \cap \overline{P}$.  

So the only issue is to find the index at which each generator $\tau_i$ of $A$ 
drops out of the filtration of higher ramification groups.  
In fact, the jumps for the $\tau_i$ are all equal.
This is because the jumps of the $\tau_i$ are the same (occurring with multiplicity $a$)
for both of the covers ${\hat{\phi}}^P$ and $\psi_\alpha$. 
To complete the proof it is thus sufficient to find the conductor 
of the cover ${\hat{\phi}}^{\alpha}$.
In particular, by Herbrand's formula, it is sufficient to show the
conductor of the $P$-Galois subcover $({\hat{\phi}}^{\alpha})^P$ equals
$\max\{s, m\sigma\}$. 

To do this, we investigate the filtration of the 
$A \times A$-Galois cover $Z \to \overline{Y}$.
Here $s$ and $m\sigma$ are the relevant upper jumps of $Z \to \overline{X}$.
Denote by $\Psi_{(\gamma, \alpha)}$ the function which takes 
the indexing on the ramification filtration of the $P \times A$-Galois
cover from the upper to lower numbering.
So the numbers $\Psi_{(\gamma, \alpha)}(s)$ and 
$\Psi_{(\gamma, \alpha)}(m\sigma)$ are the corresponding lower jumps for $Z \to \overline{X}$, 
and thus also for $Z \to \overline{Y}$ by \cite[Proposition 2]{Se:lf}.
Also the lower (and upper) jump for $Z \to \overline{Z}$ 
is $\Psi_{(\gamma, \alpha)}(m\sigma)$
and for $Z \to {\hat {Y}}$ is $\Psi_{(\gamma, \alpha)}(s)$.
Since the Galois group of $Z \to W$ is generated by 
the automorphisms $(\tau_{\gamma,i}, \tau_{\alpha,i}^{-1})$
the cover $Z \to W$ has lower jump 
$\mini\{\Psi_{(\gamma, \alpha)}(s), \Psi_{(\gamma, \alpha)}(m\sigma)\}$.

If $s \leq m\sigma$, the conductors of $Z \to W$ and of $Z \to {\hat {Y}}$ are the same;
this implies that the conductor of $W \to \overline{Y}$ equals
the conductor of ${\hat {Y}} \to \overline{Y}$, namely $\Psi_\gamma(m\sigma)$.
Thus the conductor of ${\hat{\phi}}^{\alpha}$ is $m\sigma$
and the ramification filtrations for ${\hat{\phi}}^{\alpha}$ and ${\hat{\phi}}$ 
are the same.  
Similarly, if $s > m\sigma$, the conductors of $Z \to W$ and $Z \to \overline{Z}$ are the same;
this implies that the conductor of $W \to \overline{Y}$ equals
the conductor of $\overline{Z} \to \overline{Y}$.  
So the ramification filtrations for ${\hat{\phi}}^{\alpha}$ and 
$\overline{Z} \to U$ are the same, which implies that the conductor 
of the former is also $s/m$.
\end{proof}

See Example \ref{Econd} for a computational proof of Proposition \ref{Ph1action}
in the case that $m=1$ and $P=\ZZ/p^2$.

\subsection{Connected and smooth covers}

Suppose ${\hat{\phi}}:{\hat {Y}} \to U$ is an $I$-Galois cover of normal connected germs of curves 
(wildly ramified with conductor $\sigma$ as above).  
Suppose $A \subset I^{\sigma}_{\hat{\phi}}$ satisfies the hypotheses of Lemma \ref{LhypA} and 
$\alpha \in H^\iota_A$.
In future applications \cite{Pr:deform}
and briefly in Section \ref{S3}, we need to consider questions of 
connectedness and smoothness for the $I$-Galois cover $\hat{\phi}^\alpha$.
First, we show that the $I$-Galois cover ${\hat{\phi}}^{\alpha}:W \to U$ is almost always connected.

\begin{lemma} \label{Ldisjoint}
If ${\hat{\phi}}$ is a cover of connected germs of curves, then there are only finitely many 
choices of $\alpha$ for which ${\hat{\phi}}^\alpha$ is not connected.  
In particular, ${\hat{\phi}}^\alpha$ is connected if 
${\hat{\phi}}^P$ does not dominate $\psi_\alpha$.
\end{lemma}

\begin{proof}
The first statement follows directly from the second since ${\hat{\phi}}^P$ only dominates a finite
number of $A$-Galois covers.

For the second statement, suppose ${\hat{\phi}}^P$ does not dominate $\psi_\alpha$.
The fact that both $Y$ and $V$ are invariant under the $\mu_m$-Galois action and 
the fact that $A$ is irreducible under this action imply that 
${\hat{\phi}}^P$ and $\psi_\alpha$ are linearly disjoint.
It follows that the curve $Z={\hat {Y}} \tilde{\times}_{\overline{X}} V$ is connected.
If $Z$ is connected then its quotient $W=Z^{A'}$ is connected and so 
${\hat{\phi}}^\alpha:W \to \overline{X}$ is a cover of connected germs of curves.
\end{proof}

\begin{lemma} \label{Ls=sigma}
If ${\hat{\phi}}$ is a cover of connected germs of curves and 
$\alpha$ is such that the cover ${\hat{\phi}}^\alpha$ is not connected,
then $s=m\sigma$.  
\end{lemma}

\begin{proof}
Consider the covers $\hat{\phi}:Y \to U$ and ${\hat{\phi}}^\alpha:W \to U$.
Recall that $Y$ and $W$ are two of the quotients of $Z$ by subgroups of $A \times A$.
The fact that $Y$ is connected implies that $Z$ has at most $|A|$ components.
Also $Z$ is not connected since $W$ is not connected.
The stabilizer of a component of $Z$ has size $|A|$ since $A$ is irreducible under the $\mu_m$-action.  
As a result, only one of the $A$-quotients of $Z$ can be disconnected. 
In particular, $\overline{Z}$ is connected.

The $A$-Galois covers $Y \to \overline{Y}$ and $\overline{Z} \to \overline{Y}$ are 
not linearly disjoint since their fibre product $Z \to \overline{Y}$ is disconnected. 
By the fact that $A$ is irreducible under the $\mu_m$-action, 
it follows that these two $A$-Galois covers are identical. 
Thus the $I$-Galois covers $\hat{\phi}:Y \to U$ and $\overline{Z} \to U$ are the same 
and thus have the same conductor.  The conductor of the former is $\sigma$.  
The conductor of the latter is the same as the upper jump of $V \to U$ which is $s/m$.
\end{proof}

\begin{remark} \label{Rcancel}
Note that ${\hat{\phi}}^\alpha$ could be connected even if $s =m\sigma$.
The cover ${\hat{\phi}}^{\alpha}$ is connected (and 
thus has conductor $\max\{s/m, \sigma\}$ by Proposition \ref{Ph1action})
as long as the leading term of $r_{\hat{\phi}}$ does  
not cancel the leading term of $r_\alpha +r^p-r$ for any $r \in \overline{K}$.  
For this will guarantee that 
$-\val(r_{\hat{\phi}} +r_\alpha)=\max\{-\val(r_{\hat{\phi}}), -\val(r_\alpha)\}$.
\end{remark}

A final issue that will come up is one related to smoothness.  
Fix $\Omega=\Spec(R)$ where $R$ is an equal characteristic 
complete local ring with residue field $k$ and fraction field $K$.
Consider the special fibre ${\hat{\phi}}^{\alpha}_k:W_k \to U_k$ of
the cover ${\hat{\phi}}^{\alpha}:W \to U$ of $\Omega$-curves.   
It is possible that the curve $W_k$ is singular at its closed point $w$.  
The following lemma implies that this can only occur when the ramification filtrations
of ${\hat{\phi}}^{\alpha}$ are not the same on the generic and spec.   

\begin{lemma} \label{Lsmooth}
Suppose that ${\hat{\phi}}^{\alpha}: W \to U$ is a Galois cover of normal irreducible 
germs of $\Omega$-curves with $W_k$ reduced, ${\hat{\phi}}^{\alpha}_k$ separable,
and ${\hat{\phi}}^{\alpha}$ \'etale outside $\xi$.
Let $d_k$ (resp.\ $d_K$) be the degree of the ramification divisor over
$\xi_k$ (resp.\ $\xi_K$).  
Then $d_K=d_k$ if and only if $w$ is a smooth point of $W_k$.
\end{lemma} 

\begin{proof} 
The value of $d_k$ and the question of whether $w$ is a smooth point of $W_k$ do not depend
on $R$.  The value of $d_K$ depends only on the ramification at the generic point
of $\xi_K$.
For this reason, it is sufficient to restrict to the case that $R=k[[t]]$.
The proof then follows from Kato's formula \cite{Ka}
which states (under the conditions above) that
$\mu_{w}-1=n(\mu_{\xi_k}-1)+d_K-d_k^{wild}$.
Here $n=\deg({\hat{\phi}}^{\alpha})$, $\mu_{\xi_k}=0$ since ${\xi_k}$ is a smooth point of $U_k$
and $\mu_{w}$ is the invariant measuring the singularity of $W_k$ at $w$. 
Namely, let $\displaystyle \pi_w:\tilde{W_k}
\to W_k$ be the normalization of $W_k$.
Define $\mu_{w}=2\delta_{w}-m_{w}+1$
where $\delta_{w}=\dime_k(\tilde{\CO}_w/\CO_w)$ and 
$m_w=|\pi_w^{-1}(w)|$.
We see that $\mu_w=0$ if and only if $w$ is smooth on $W_k$. 

Write $\displaystyle \tilde{\phi}_k:\tilde{W_k} \to U_k$
and $\displaystyle \tilde{\phi}^{-1}_k(\xi_k)=\{w_i| 1 \leq i \leq m_w\}$.
Let $\delta_i$ and $e_i$ be respectively the discriminant ideal 
and the ramification index of the extension $k[[u]] \to \hat{\CO}_{\tilde{W}_k,w_i}$
on the special fibre.
Let $d_i=v(\delta_i)$ and let $d_i^{wild}=d_i-e_i+1$. 
Let $d_k=\sum_i d_i$ and let $d_k^{wild}=\sum_i d_i^{wild}$.

By Kato's formula, 
$\mu_w=1-n+d_K-d_k^{wild}=1-n+ d_K -m_w(v(\delta_i) -e_i +1)$.
Thus $\mu_w= 1-n+d_K -d_k +n -m_w=1-m_w +d_K -d_k$.
Note that $d_K \geq d_k$.
So the condition that $w$ is smooth implies 
$\mu_w=0$, $m_w=1$, and $d_K=d_k$.
Conversely, if $d_K=d_k$ then $\mu_w=1-m_w$, 
but $m_w \geq 1$ so $\mu_w=0$ and $w$ is smooth on $W_k$.
\end{proof}

\section{Increasing the Conductor of a Cover.} \label{S3}

In this section, we consider the question of increasing the conductor 
of a cover of curves at a wild branch point.  
Suppose $X$ is a smooth projective irreducible $k$-curve and 
$B$ is a non-empty finite set of points of $X$.
Suppose $G$ is a finite quotient of $\pi_1(X-B)$. 
When $|B|$ is nonempty, these groups have been classified by 
Raynaud and Harbater (\cite{Ra:ab} \cite{Ha:ab}) in their proof of Abhyankar's Conjecture.
Namely, $G$ is a finite quotient of $\pi_1(X-B)$ if and only if the number of generators 
of $G/p(G)$ is at most $2g_X+|B|-1$.  Here $p(G)$ denotes the characteristic subgroup of $G$
generated by the group elements of $p$-power order. 

Suppose $\phi:Y \to X$ is a $G$-Galois cover of smooth projective irreducible 
$k$-curves branched only at $B$.
Suppose $\phi$ wildly ramified over $b \in X$ and 
let $I$ be the inertia group at a ramification point $\eta$ above $b$. 
Let $R$ be an equal characteristic complete discrete valuation ring with 
residue field $k$ and let $K=\Frac(R)$.
In this section, we will ``deform'' $\phi$ to a cover $\phi_R:Y_R \to X \times_k R$
so that the conductor of $\phi_R$ at this wild branch point increases on the generic fibre.
By Lemma \ref{Lsmooth}, it is not possible to do this without introducing a singularity.
In other words, the special fibre $\phi_k$ of $\phi_R$ will be singular and 
$\phi$ will be isomorphic to the normalization of $\phi_k$ away from $b$.  

More precisely, let $\hat{\phi}:\hat{Y} \to U_k$ be the germ of $\phi$ at $\eta$.  
Here $U_k \simeq \Spec(k[[u]])$.  
Let $U_R=\Spec(k[[t]][[u]])$.  Let $\xi$ (resp. $\xi_R$) 
be the closed point of $U_k$ (resp.\ $U_R$) given by the equation $u=0$.
For lack of better terminology, a {\it singular deformation} of $\hat{\phi}$
is an $I$-Galois cover $\hat{\phi}_R: \Y \to U_R$ of normal irreducible germs of $R$-curves,  
whose branch locus consists of only the $R$-point $\xi_R$, such that 
the normalization of the special fibre of $\hat{\phi}_R$ is isomorphic to 
$\hat{\phi}$ away from $\xi$.  
(Note that the {\it special fibre} $\hat{\phi}_k$ of $\hat{\phi}_R$ is the restriction of 
$\hat{\phi}_R$ over $U_k$ (when $t=0$).
The {\it generic fibre} $\hat{\phi}_K$ of $\hat{\phi}_R$ is the cover of 
$U_R \times_R K=\Spec(k[[t]][[u]][t^{-1}])$.

Lemma \ref{Lsmooth} shows that a deformation is singular if and only if the ramification 
data on the generic and special fibres are not the same.

The following proposition shows that there is
always a singular deformation of $\hat{\phi}$ with larger conductor on the generic fibre.

\begin{proposition}  \label{Pupdef}
Suppose $\hat{\phi}: \hat{Y} \to U_k$ is an $I$-Galois cover of normal 
connected germs of curves with conductor $\sigma$.   
Suppose $A \subset I^\sigma$ satisfies the hypotheses of Lemma \ref{LhypA}.
Suppose $s \in \NN$ is such that $p \nmid s$, $s > m\sigma$, and  
$s \equiv s_\iota \bmod m$ (with $s_\iota$ as in Lemma \ref{Lcong}). 
Then there exists a singular deformation $\hat{\phi}_R: \Y \to U_R$ 
whose generic fibre $\hat{\phi}_K:\Y_K \to U_K$
has inertia $I$ and conductor $s/m$.
In addition, $I_{\phi_K}^c=I_{\phi}^c$ 
for $0 \leq c \leq \sigma$, $I_{\phi_K}^c=A$ 
for $\sigma < c \leq s/m$, and $I_{\phi_K}^c=0$ for $c > s/m$.
\end{proposition}

\begin{proof} 
The $A$-Galois subcover $\hat{Y} \to \overline{Y}$ 
is determined by $\kappa$ and by the equation $y^p-y=r_\phi$
of its $\langle \tau_1 \rangle$-Galois quotient
(where $r_\phi \in \overline{K}$).  
Consider the singular deformation of $\hat{\phi}$ whose $\overline{P}$-Galois quotient 
is constant and whose $\langle \tau_1 \rangle$-Galois subquotient 
is given generically by the following equation: $y^p-y=r_\phi + tx^{-s}$.
The curve $\Y$ is singular only above the point $(u,t)=(0,0)$.  
The normalization of the special fibre agrees with $\hat{\phi}$ away from $u=0$.
The cover $\hat{\phi}_R$ is branched only at $\xi_R$ since $u=0$ 
is the only pole of the function $r_\phi + tx^{-s}$.
When $t \not = 0$, by Proposition \ref{Ph1action}, 
the conductor of $\hat{\phi}_K$ is $s$ 
and the statement on ramification filtrations is true.
\end{proof}

Proposition \ref{Pupdef} shows that we can increase the conductor from 
$\sigma$ to $s/m$.  By Herbrand's formula, it follows that the last lower jump 
increases from $j_e$ to $j_e +p^{e-a}m(s/m-\sigma)$.  One can check that the latter number
is always an integer which is not divisible by $p$.

The next theorem uses Proposition \ref{Pupdef} and formal patching
to (singularly) deform a cover $\phi:Y \to X$ of projective curves to a family of
covers $\phi_R$ of $X$ so that the conductor increases at the chosen branch point.  
This family can be defined over a 
variety $\Theta$ of finite type over $k$.  We then specialize to
a fibre of the family over another $k$-point of $\Theta$ to get a cover
$\phi'$ with larger conductor.

\begin{theorem} \label{2patch1}
Let $\phi: Y \to X$ be a $G$-Galois cover of smooth projective irreducible
curves with branch locus $B$.
Suppose $\phi$ is wildly ramified with inertia group $I \simeq P \rtimes_\iota \mu_m$ 
and conductor $\sigma$ above some point $b \in B$.  
Suppose $A \subset I^\sigma$ satisfies the hypotheses of Lemma \ref{LhypA}.
Let $s_\iota$ be as defined in Lemma \ref{Lcong}(ii).   
Suppose $s \in \NN$ such that $p \nmid s$, $s > m\sigma$, and $s \equiv s_\iota \bmod m$. 
Then there exist
$G$-Galois covers $\phi_R: Y_R \to X \times_k R$ and $\phi': Y' \to X$ such that:
\begin{enumerate}
\item The curves $Y_R$ and $Y'$ are irreducible and $Y_K$ and $Y'$ are smooth
and connected. 

\item After normalization, the special fibre $\phi_k$ of $\phi_R$
is isomorphic to $\phi$ away from $b$. 

\item The branch locus of the cover $\phi_R$ (resp.\ $\phi'$) 
consists exactly of the $R$-points $\xi_R=\xi \times_k R$ 
(resp.\ the $k$-points $\xi$) for $\xi \in B$.

\item For $\xi \in B$, $\xi \not = b$, 
the ramification behavior for $\phi_R$ at $\xi_R$ 
(resp.\ $\phi'$ at $\xi$) is identical to that of $\phi$ at $\xi$.

\item At the $K$-point $b_K$ (resp.\ at $b$),
the cover $\phi_K$ (respectively $\phi'$) 
has inertia $I$ and conductor $s/m$ and in addition
$I_{\phi_K}^c=I_{\phi}^c$ for $0 \leq c \leq \sigma$, $I_{\phi_K}^c=A$ 
for $\sigma < c \leq s/m$, and $I_{\phi_K}^c=0$ for $c > s/m$.

\item The genus of $Y'$ and of $Y_K$ is 
$g'_Y= g_Y + |G|(s/m-\sigma)(1-1/p^a)/2$ where $a$ is such that $|A|=p^a$. 
\end{enumerate}  
\end{theorem}

\begin{proof}
Let $\eta \in \phi^{-1}(b)$ and
consider the $I$-Galois cover $\hat{\phi}: \hat{Y}_{\eta} \to \hat{X}_{b}$.
Applying Proposition \ref{Pupdef} to $\hat{\phi}$, 
there exists a singular deformation $\hat{\phi}_R: \hat{Y}_R \to \hat{X}_R$ of 
$\hat{\phi}$ with the desired properties.  In particular, $\hat{\phi}_K$ has 
inertia $I$ and conductor $s$ over $b_K$.  
Consider the disconnected $G$-Galois cover $\Ind_I^G(\hat{\phi}_R)$.

The covers $\phi_k$ and $\Ind_I^G(\hat{\phi}_R)$ and the isomorphism
between their overlap 
constitute a relative $G$-Galois thickening problem.  
The (unique) solution to this thickening problem \cite[Theorem 4]{HS} 
yields the $G$-Galois cover $\phi_R:Y_R \to X \times_k R$.  
Recall that the cover $\phi_R$ is isomorphic to
$\Ind_I^G(\hat{\phi}_R)$ over $\hat{X}_{R,b}$.  Also, 
$\phi_R$ is isomorphic to the 
trivial deformation $\phi_{tr}:Y_{tr} \to X_{tr}$ of $\phi$ away from $b$.
Thus $Y_R$ is irreducible since $Y$ is irreducible and 
$Y_K$ is smooth since $Y_{tr,K}$ and $\hat{Y}_K$ are smooth.  

The data for the cover $\phi_R$ is contained in a 
subring $\Theta \subset R$ of finite 
type over $k$, with $\Theta \not = k$ since the family is non-constant.  
Since $k$ is algebraically closed, there exist infinitely many $k$-points of $\Spec(\Theta)$.  
The cover $\phi_R$ descends to a cover of $\Theta$-curves.
The closure $L$ of the locus of $k$-points $\theta$ of $\Spec(\Theta)$ over 
which the fibre $\phi_\theta$ is not a $G$-Galois cover of 
smooth connected curves is closed, \cite[Proposition 9.29]{FJ}.  
Furthermore, $L \not = \Spec(\Theta)$ since $Y_K$ is smooth and irreducible.
Let $\phi':Y' \to X$ be the fibre over a $k$-point not in $L$.  
Note that $Y'$ is smooth and irreducible by definition.

Properties 2-5 follow immediately from the compatibility of $\phi_R$ with 
$\Ind_I^G(\hat{\phi}_R)$ over $\hat{X}_{R,b}$ and with the trivial deformation 
$\phi_{tr}:Y_{tr} \to X_{tr}$ away from $b$.

The genus of $Y'$ and of $Y_K$ 
increases because of the extra contribution to the Riemann-Hurwitz formula.  
In particular, there are $|G|/m|P|$ ramification points above $b_K$. 
If $|P|=p^e$ then by Herbrand's formula, each ramification point has 
$(s/m-\sigma)mp^{e-a}$ extra non-trivial ramification groups of order
$p^a$ in the lower numbering.  
Thus the degree of the ramification divisor over $b_K$ increases by $|G|(s/m-\sigma)(1-1/p^a)/2$. 
\end{proof}

One can say more when $X \simeq \PP_k$ and $B=\{\infty\}$.
By Abhyankar's Conjecture, the non-trivial quasi-$p$ groups are exactly the ones so that 
$G$ is a finite quotient of $\pi_1(\Aa_k)$ and $p$ divides $|G|$. 

\begin{corollary} \label{Ca1jump}
If $G \not = 0$ is a quasi-$p$ group and $\sigma \in \NN$ ($p \nmid \sigma$)
is sufficiently large, 
then there exists a $G$-Galois cover $\phi:Y \to \PP_k$ branched at only one 
point with conductor $\sigma$.   
\end{corollary}

The author obtained a similar result in \cite{Pr:cond2} under the 
restriction that the Sylow $p$-subgroup of $G$ has order $p$.

\begin{proof}
By \cite{Ra:ab} and \cite{Ha:ab}, there exists a $G$-Galois cover $\phi:Y \to \PP_k$
branched only at $\infty$ whose inertia groups are the Sylow $p$-subgroups of $G$.
The result is then automatic from Theorem \ref{2patch1}.
\end{proof}

As another corollary, we show that curves $Y$ of arbitrarily high genus occur 
for $G$-Galois covers $\phi:Y \to X$ branched at $B$ as long as $p$ divides $|G|$.

\begin{corollary} \label{CbigG}
Suppose $X$ is a smooth projective irreducible $k$-curve and $B \subset X$ is a non-empty
finite set of points.  Suppose $G$ is a finite quotient of $\pi_1(X-B)$ such that
$p$ divides $|G|$.  
Let ${\mathfrak p} \subset \NN$ be the set of genera $g$ for which  
there exists a $G$-Galois cover $\phi:Y' \to X$ branched only at $B$ with $\gen(Y')=g$.
Then the set ${\mathfrak p}$ contains an arithmetic progression whose increment 
depends only on $G$ and $p$.  
\end{corollary}

\begin{proof}
First we show that the hypotheses on $G$ guarantee the existence of a $G$-Galois cover 
$\phi:Y \to X$ branched only at $B$ with wild ramification at some point.  
Let $p(G)$ be the normal subgroup of $G$ generated by all elements of $p$-power order.
Let $S \subset p(G)$ be a Sylow $p$-subgroup of $G$.  
Consider the natural morphism $\pi:G \to G/p(G)$.
By \cite[Lemma 2.4]{H:emAb}, there exists $F \subset G$ which is prime-to-$p$ and 
normalizes $S$ so that $\pi(F)=G/p(G)$.
Let $g_1, \ldots, g_r$ be a minimal set of generators for $G/p(G)$.
After possibly replacing $F$ with the subgroup generated by elements $h_1, \ldots, h_r$ 
where $\pi(h_i)=g_i$, we see that $F$ can be generated by $r$ elements.
By Abhyankar's Conjecture \cite{Ha:ab}, $r \leq 2g_X+|B|-1$ and there exists an $F$-Galois cover
$Y_\circ \to X$ branched only over $B$.    
Note that $F$ and $p(G)$ generate $G$.
As a result, the $F$-Galois cover $Y_\circ \to X$ and $p(G)$ satisfy all the hypotheses of 
\cite[Theorems 2.1 and 4.1]{H:emAb}.  Let $I_1$ be the inertia
group of $Y_\circ \to X$ above a chosen point $b \in B$.
These theorems allow one to modify the cover $Y_\circ \to X$ to get a new 
$G$-Galois cover $\phi:Y \to X$ branched only at $B$ so that the inertia above 
$b$ is $I=I_1S$.  Since $p$ divides $|G|$, it follows that $S$ is non-trivial 
and so $\phi$ is wildly ramified above $b$.  

Let $g_Y$ be the genus of $Y$.  The inertia group $I$ is of the form $P \rtimes_\iota \mu_m$
for some $P \subset S$.  Suppose $|P|=p^e$.  Let $\sigma$ be the conductor of $\phi$ above $b$.
Let $n$ be the fixed integer $|G|(1-1/p^a)\sigma/2$.
Let $s$ be such that $p \nmid s$, $s > m\sigma$, and $s \equiv s_\iota \bmod m$. 
By Theorem \ref{2patch1}, it is possible to deform $\phi$ to produce another
curve $\phi':Y'\to X$ branched only at $B$ with larger genus $g_Y'=g_Y -n + |G|(1-1/p^a)s/2m$.  
The set of $g_Y'$ realized in this way clearly contains an arithmetic progression. 

By the congruence condition on $s$, one can increase $s$ only 
by a multiple $s'm$ of $m$, which causes the genus to increase by $|G|(1-1/p^a)s'/2$.  
But one has to remove $(1/p)$th of these values since
the integer $s+s'm$ will be divisible by $p$ exactly $(1/p)$th of the time.
So the set ${\mathfrak p}$ contains $p-1$ arithmetic progressions with increment $p|G|(1-1/p^a)/2$. 
\end{proof}

\begin{remark}
The proportion of the set ${\mathfrak p}$ in $\NN$ is at least $2(p^a-p^{a-1})/|G|(p^a-1)$
for some $a \geq 1$ 
such that the center of $S$ contains a subgroup $A \simeq (\ZZ/p)^a$. 
This lower bound is approximately $2/|G|$ for large $p$.  
It is realized when $X \simeq \PP_k$, $B=\{\infty\}$ and $G=\ZZ/p$. 
To see this, note that $\ZZ/p$-covers of the affine line correspond to Artin-Schreier equations 
$y^p-y=f(x)$ where the degree $j$ of $f(x)$ is prime-to-$p$.
Such a curve has genus $(p-1)(j-1)/2$.  
It follows that the proportion of genera which occur in this case is exactly $2/p$.
But in general, we expect that this lower bound is not optimal.  
For example, we expect that this proportion equals $2/p$ whenever $G$ is an abelian $p$-group.
\end{remark}

The following corollary shows that the 
the structure of Hurwitz spaces for wildly ramified covers will be vastly different
from those of tamely ramified covers.  

\begin{corollary} \label{Churwitz}
For any smooth projective irreducible $k$-curve $X$
and any non-empty finite set of points $B \subset X$
and any finite quotient of $\pi_1(X-B)$ so that $p$ divides $|G|$,   
a Hurwitz space for $G$-Galois covers $\phi:Y \to X$ branched at $B$ will have infinitely 
many components.
\end{corollary}

\begin{proof}
The proof is immediate from Corollary \ref{CbigG} since 
two covers with different genus cannot correspond to points in the same component 
of a Hurwitz space.
\end{proof}

\section{Example: Inertia $\ZZ/p^e$.} \label{Sex}

In this section, we study $\ZZ/p^e$-Galois covers using class field theory.  
We give an example of Proposition \ref{Ph1action} and of singular deformations.

\begin{definition}
A sequence $\sigma_1, \ldots, \sigma_e$ is {\it $p^e$-admissible}
if $\sigma_i \in \NN^+$, $p \nmid \sigma_1$ and for $1 \leq i \leq e-1$, 
either $\sigma_{i+1} =p \sigma_i$ or $\sigma_{i+1} > p \sigma_i$ and 
$p \nmid \sigma_{i+1}$. 
\end{definition}

For convenience, we include the proof of the 
following classical lemma, \cite{Schmid}.  

\begin{lemma} \label{p^ejumps}
If $\hat{\phi}:\hat{Y} \to U_k$ is a $\ZZ/p^e$-Galois cover, then the upper jumps 
$\sigma_1, \ldots, \sigma_e$ of its ramification filtration are $p^e$-admissible.    
For any $p^e$-admissible sequence $\Sigma$, there exists a $\ZZ/p^e$-Galois cover 
$\hat{\phi}:\hat{Y} \to U_k$ with upper jumps at the indices in $\Sigma$. 
\end{lemma}

\begin{proof} By the Hasse-Arf Theorem, $\sigma_i \in \NN^+$.
Since $\sigma_1=j_1$, we see that $p \nmid \sigma_1$.
Suppose $k((u)) \hookrightarrow L$ is the $p^e$-Galois field extension 
corresponding to $\hat{\phi}$.  
Note that for $\sigma_i < n \leq \sigma_{i+1}$, 
the $n$th ramification group $I^n$ in the upper numbering 
equals $p^i\ZZ/p^e$.
Denote by $U^n$ the unit group $(1+u^nk[[u]]) \subset k[[u]]^*$.
By local class field theory, there is a reciprocity isomorphism 
$\omega:k((u))^*/NL^* \to \ZZ/p^e$ so that the image of $U^n$ under $\omega$
equals $I^n$, \cite[Chapter XV]{Se:lf}.

First we show that $\sigma_{i+1} \geq p\sigma_i$.
There is some $1+t^{\sigma_i}h \in U^{\sigma_i}$ whose image under $\omega$
generates $p^{i-1}\ZZ/p^e$.  Thus the image of 
$(1+t^{\sigma_i}h)^p \equiv 1+t^{p\sigma_i}h^p \in U^{p\sigma_i}$
generates $pp^{i-1}\ZZ/p^e$.  Thus $p^i\ZZ/p^e \subset I^{p\sigma_i}$
which implies $\sigma_{i+1} \geq p\sigma_i$. 

Next suppose $p|\sigma_{i+1}$ and write $\sigma_{i+1}=p(\sigma_i +c)$
for some $c \in \NN$.  The image of $1+t^{p(\sigma_i +c)} \in U^{\sigma_{i+1}}$
must generate $p^i\ZZ/p^e$ so the image of $1+t^{\sigma_i +c}$ must 
generate $p^{i-1}\ZZ/p^e$.  Thus $1+t^{\sigma_i+c} \in U^{\sigma_i}$
which implies $c=0$.     

Finally, for any $p^e$-admissible sequence $\{\sigma_1, \ldots, \sigma_e\}$
we define a homomorphism $U^1 \to \ZZ/p^e$ as follows:
if $p \nmid c$ and $\sigma_i < c \leq \sigma_{i+1}$, 
then $1+t^{p^rc} \mapsto p^rp^i \in \ZZ/p^e$.  Since $U^1$ is isomorphic 
to the abelianization of the fundamental group of $k((u))$, it follows that
there exists a $\ZZ/p^e$-Galois cover $\hat{\phi}:\hat{Y} \to U_k$.
The upper jumps of $\hat{\phi}$ are the given sequence since
$I^n \not = I^{n+1}$ if and only if $n=\sigma_i$ for some $i$.
\end{proof}

We define a partial ordering on the set of {\it $p^e$-admissible}
sequences as follows.
 
\begin{definition}
Suppose $\Sigma$ and $\Sigma'$ are two $p^e$-admissible sequences
given by $\sigma_1, \ldots, \sigma_e$ and $\sigma'_1, \ldots, \sigma'_e$
respectively.
Then $\Sigma' > \Sigma$ if $\sigma'_i > \sigma_i$ for $1 \leq i \leq e$.
\end{definition}

The $p^e$-admissible sequence which is smaller than all others is  
$\{1,p,p^2, \ldots, p^{e-1}\}$.  

In the following example, we give a computational proof of Proposition \ref{Ph1action}
in the case that $m=1$ and $P=\ZZ/p^2$.  We note that, even in this simple case, the equations
are quite complicated.

\begin{example} \label{Econd}
Suppose $P=\ZZ/p^2$ and $A=\ZZ/p$.  
Choose $j \in \NN^+$ with $p \nmid j$.
By Lemma \ref{p^ejumps}, there exists a $P$-Galois cover
${\hat{\phi}}:{\hat {Y}} \to \overline{X}$ 
of germs of $k$-curves with upper jumps $\sigma_1=j$ and $\sigma_2=pj$.  
Furthermore, no conductor smaller than $pj$ can occur for a $P$-Galois cover with $\sigma_1=j$.
By Herbrand's formula, 
the lower jumps of ${\hat{\phi}}$ are $j_1=j$ and $j_2=(p^2-p+1)j$.

After some changes of coordinates, the equation for the $\ZZ/p$-Galois 
quotient $\overline{\phi}:\overline{Y}\to \overline{X}$ can be given generically 
by $y_1^p-y_1=x^{-j}$.
Note that there is a natural valuation $\val$ on the fraction field
$\overline{K}=k((x))[y_1]/(y_1^p-y_1-x^{-j})$ of $\overline{Y}$ and 
$\val(x)=p$ and $\val(y_1)=-j$. 
The equation for the $\ZZ/p$-Galois subcover ${\phi^A}$ is $y_2^p-y_2=f(y_1,x)$ 
for some $f(y_1, x) \in \overline{K}$ with valuation $-j_2$.   

The $\ZZ/p$-Galois cover $\psi_\alpha$ is given by an equation 
$v^p-v=g(x)$ for some $g(x) \in k((x))$.
Let $s$ be the degree of $x^{-1}$ in $g(x)$.
By Lemma \ref{Lfibre}, the cover ${\hat{\phi}}^{\alpha}$
is determined by its $\ZZ/p$-Galois subcover $W \to \overline{Y}$
which is given by the equation $w^p-w=f(y_1,x) +g(x)$.
By Lemma \ref{Ltransitive1}, 
any $P$-Galois cover dominating $\overline{\phi}$ is of this form for some 
$g(x) \in k((x))$.  

We now give an explicit proof of Proposition \ref{Ph1action} in this case.
Namely, we show that the upper jumps of ${\hat{\phi}}^{\alpha}$ are $\sigma_1$ and 
$\max\{s, \sigma_2\}$.  The first upper jump is $\sigma_1$ since 
${\hat{\phi}}^{\alpha}$ dominates $\overline{\phi}$.
Let $J$ be the last lower jump of ${\hat{\phi}}^{\alpha}$. 
Recall that $J$ is the prime-to-$p$ valuation of $f(y_1,x) +g(x)$.
Note that $x^{-s}=(y_1^p-y_1)^{s/j}=y_1^{ps/j}(1-y_1^{1-p})^{s/j}$.
So $$x^{-s}=y_1^{ps/j}(1- y_1^{-(p-1)}s/j + \ldots)=
y_1^{ps/j} - y_1^{ps/j -(p-1)}s/j + \ldots.$$
We can modify the equation for $W \to \overline{Y}$
by adding $-y_1^{ps/j} + y_1^{s/j}$ without 
changing its isomorphism class as a $\ZZ/p$-Galois cover or its jump.
The next term in the above equation indicates that
the lower jump $J$ of the cover $w^p-w=f(y_1,x) +g(x)$ equals
$\max\{-\val(y_1^{ps/j -(p-1)}), -\val(f(y_1,x))\}$.  (It cannot be smaller since 
the valuation of $f(y_1,x)$ is minimal.)
So $J=\max\{ps-j(p-1), j_2\}$. 
By Herbrand's formula, the conductor of ${\hat{\phi}}^{\alpha}$ equals $\max\{s, \sigma_2\}$.
\end{example}

\begin{proposition}  \label{Pupdef2}
Suppose there exists a $\ZZ/p^e$-Galois cover $\hat{\phi}: \hat{Y} \to U_k$ of normal 
connected germs of curves whose upper jumps in the ramification filtration
are at the $p^e$-admissible sequence $\Sigma$.  Suppose $\Sigma' > \Sigma$
is a $p^e$-admissible sequence.     
Then there exists a singular deformation $\hat{\phi}_\Omega: \Y \to U_{\Omega}$ 
whose generic fibre has ramification filtration $\Sigma'$.
\end{proposition}

\begin{proof} 
The proof is by induction on $e$.  If $e=1$, the proof follows by
Proposition \ref{Pupdef}.  If $e > 1$, choose 
$A \subset I^{\sigma}$ to be a subgroup of order $p$.
The $P/A=\overline{P}$-Galois quotient $\overline{\phi}$ has upper jumps
$\sigma_1, \ldots \sigma_{e-1}$.  By the inductive hypothesis, 
there exists a deformation $\overline{\phi}_\Omega:\overline{Y}_\Omega \to U_{\Omega}$
whose generic fibre is a $\ZZ/p^{e-1}$-Galois cover of normal connected curves
whose branch locus consists of only the $K$-point $\xi_K=\xi_\Omega \times_\Omega K$
over which it has upper jumps $\sigma'_1, \ldots, \sigma'_{e-1}$.

By \cite[X, Theorem 5.1]{AGV}, there exists a $P$-Galois cover $\hat{\phi}'_\Omega$
dominating $\overline{\phi}_\Omega$.  
Choose $\hat{\phi}'_\Omega$ to have minimal conductor $s'$ among all such covers dominating
$\overline{\phi}_\Omega$. 
By Lemma \ref{p^ejumps}, $s'=p\sigma'_{e-1}$. 
The restriction of $\hat{\phi}'_\Omega$ to $U_k$ has conductor at most $s'$. 
This restriction differs from  $\phi$ by an element 
$\alpha \in \Hom(\pi_1(U_\Omega - \xi_\Omega), A)$.
By Proposition \ref{Ph1action}, $\alpha$ has conductor at most $\max\{s',\sigma_e\}$.

Let $\hat{\phi}_e$ be the cover $\hat{\phi}'_\Omega$ modified by $\alpha$.
By Proposition \ref{Ph1action} and by minimality of $s'$,
the conductor $\sigma$ of $\hat{\phi}_e$ satisfies
$s' \leq \sigma \leq \max\{s',\sigma_e\}$.
Since $\Sigma'>\Sigma$ is $p^e$-admissible, $\sigma'_e \geq \sigma_e$
and $\sigma'_e \geq p \sigma'_{e-1}=s'$ (with $p \nmid \sigma'_e$ if equality 
does not hold).  We apply Proposition \ref{Pupdef} with $s=\sigma'_e$
to increase the conductor.  The conclusion is that 
there exists a deformation $\hat{\phi}_\Omega: \Y \to U_{\Omega}$ 
whose generic fibre has inertia $\ZZ/p^e$ and conductor $\sigma'_e$.
\end{proof} 

\bibliographystyle{abbrv} 
\bibliography{paper}

\noindent
Rachel J. Pries\\
Department of Mathematics\\
Colorado State University\\
Fort Collins, CO 80523\\
pries@math.colostate.edu

\end{document}